\documentclass[twocolumn,10pt]{asme2emod}

\usepackage{graphicx}
\usepackage{float}
\usepackage{amsmath}
\usepackage{amsfonts}
\usepackage{amssymb}
\usepackage{wrapfig}
\usepackage[thinspace, thinqspace,amssymb]{SIunits}
\usepackage{subfigure}
\usepackage{placeins}
\usepackage{acronym}
\usepackage[hidelinks]{hyperref}
\usepackage{pgfplots}
\usepackage{tikz}
\usepackage{tikz}
\usetikzlibrary{shapes,arrows}
\usetikzlibrary{positioning}
\title{Adjoint Method to Calculate the Shape Gradients of Failure Probabilities for Turbomachinery Components}

\author{Hanno Gottschalk\\{\tensfb Mohamed Saadi}\\{\tensfb Onur Tanil Doganay}\\{\tensfb  Kathrin Klamroth}
    \affiliation{School of Mathematics and Natural Science\\
    University of Wuppertal,\\ Gaussstr 20, 42119 Wuppertal,\\
    Germany\\
    Email: hanno.gottschalk@uni-wuppertal.de\\
    saadi@uni-wuppertal.de\\
    doganay@uni-wuppertal.de\\
    klamroth@math.uni-wuppertal.de
    }\\

}

\author{  {\tensfb Sebastian Schmitz}
    \affiliation{Gas Turbine Department of Materials and Technology\\
		Siemens AG\\
       Berlin, 10553\\
    Germany\\
    Email: schmitz.sebastian@siemens.com\vspace{5mm}
    }\\

}

\begin{document}

\maketitle    

%%%%%%%%%%%%%%%%%%%%%%%%%%%%%%%%%%%%%%%%%%%%%%%%%%%%%%%%%%%%%%%%%%%%%%
\begin{abstract}
{\it In the optimization of turbomachinery components, shape sensitivities for fluid dynamical objective functions have been used for a long time. As peak stress is not a differential functional of the shape, such highly efficient procedures so far have been missing for objective functionals that stem from mechanical integrity. This changes, if deterministic lifing criteria are replaced by probabilistic criteria, which have been introduced recently to the field of low cycle fatigue (LCF).
 
Here we present a finite element (FEA) based first discretize, then adjoin approach to the calculation of shape gradients (sensitivities) for the failure probability with regard to probabilistic LCF and apply it to simple and complex geometries, as e.g. a blisk geometry.

We review the computation of failure probabilities with a FEA postprocessor and sketch the computation of the relevant quantities for the adjoint method. We demonstrate high accuracy and computational efficiency of the adjoint method  compared to finite difference schemes. We discuss implementation details for rotating components with cyclic boundary conditions. Finally, we shortly comment on future development steps and on potential applications in multi criteria optimization.  
}

\end{abstract}

%%%%%%%%%%%%%%%%%%%%%%%%%%%%%%%%%%%%%%%%%%%%%%%%%%%%%%%%%%%%%%%%%%%%%%
\begin{nomenclature}
\entry{$LCF$}{Low cycle fatigue}
\entry{$SBO$}{Surrogate based optimization}
\entry{$CFD$}{Computational fluid dynamics}
\entry{$\Omega$}{3D component domain}
\entry{$\partial\Omega$}{2D domain surface}
\entry{$\lambda,\mu$}{Lam\'e's constants for linear elasticity}
\entry{$FE$, $FEA$}{finite element (analysis)}
\entry{$X$}{finite element node set}
\entry{$DoF$}{Degrees of freedom}
\entry{$N$,$M$,$q$}{Number of global Nodes/elements/local DoF}
\entry{$u,U$}{Displacement field}
\entry{$\nabla u$}{Jacobi matrix of $u$}
\entry{$\theta_j(x)$, $\hat\theta_j(\hat x)$}{finite element shape functions}
\entry{$B$}{finite element stiffness matrix}
\entry{$f$,$g$,$F$}{finite element volume and surface load densities}
\entry{$\sigma$}{Stress tensor (field)}
\entry{$\varepsilon$}{strain tensor (field)}
\entry{$\sigma_{\rm vM}$}{Equivalent (elastic) stress}
\entry{$CMB$}{Coffin-Manson-Basquin (model)}
\entry{$\sigma_f$}{Fatigue strength coefficient}
\entry{$b$}{Fatigue strength exponent}
\entry{$\epsilon_f$}{Fatigue ductility coefficient}
\entry{$c$}{Fatigue ductility exponent}
\entry{$E$}{Cyclic Young's modulus}
\entry{$N_i$, $n$}{Load cycles until crack initiation and cycle count}
\entry{$PoF$}{Probability of failure}
\entry{$m$}{Weibull shape parameter}
\entry{$\eta$}{Weibull scale parameter}
\entry{$\xi_j$,$\hat\xi_j$, $\hat\omega_j$}{quadrature points and weights (surface and volume)}
\entry{$\Lambda$}{adjoint state}
\entry{$\mathcal{L}$}{Lagrange function}
\entry{$\omega$}{angular velocity}
\entry{$\varrho$}{density}
\entry{$\chi^*$}{Normalized gradient of equivalent elastic stress}
\end{nomenclature}

%%%%%%%%%%%%%%%%%%%%%%%%%%%%%%%%%%%%%%%%%%%%%%%%%%%%%%%%%%%%%%%%%%%%%%
\section{Introduction}
For more than one decade, algorithmic optimization of gas turbine components has been consistently applied to improve the efficiency and reliability, see  \cite{amtsfeld,kim, buske, backhaus,xu} to name just a few examples. Two major trends can be identified in the literature: Surrogate based optimization (SBO), see e.g. \cite{amtsfeld,kim,forrester}, and gradient based methods using the adjoint approach \cite{xu}. In some works \cite{forrester,backhaus}, gradient enhanced kriging \cite{debar} has been applied for a combination of both trends.

Analyzing the technical state of the art, a bias between theoretical and practical preferences in the selection of efficient optimization algorithms becomes obvious. While from a theoretical standpoint, the adjoint mehthod should be preferred due to guaranteed convergence to (local) minima along with error estimates \cite{nocedal,troeltsch} and avoidance of the \emph{curse of dimension}, practitioners mostly prefer SBO, see e.g. \cite{forrester}. This is even more surprising, as efficient adjoint fluid dynamics codes are available, see e.g. \cite{giles1,giles2,wever,frey}.          

The authors propose three main reasons to explain the above situation. First, adjoint CFD solvers are somewhat sensitive to CFD settings, as e.g. the size of residuals in the iterative solvers of the flow field. Consequently, shape gradients (also called shape sensitivities) calculated with the adjoint method are somewhat rough on some parts of the geometry, as e.g. trailing edges of vanes and blades. Second, adjoint CFD solvers have to be adjusted whenever the baseline solver is improved, which can be a time and resource intensive process. Third, up to now there has been no clear recipe how to combine the adjoint method from CFD with requirements of structural integrity. Therefore, the look ahead to multi-physics, multi-objective optimization seems to be better understood in the case of SBO, where several objectives can be treated on the same footing. 

While all three reasons are valid, it is the third reason that can not be overcome by software technologies, like gradient redefinition in the first case, see e.g. \cite{schulz}, or automated differentiation (AD) for the second case \cite{nocedal}. Let us consider, e.g., the deterministic design life for turbine blades and vanes with regard to Low Cycle Fatigue (LCF). As the deterministic design life is calculated at the point of the highest loading, the safe number of cycles is calculated for all surface points of the component and then is minimized. The operation of minimization, however, is not differentiable. This is more than just mathematical sophistry, as the location of the point of    highest loading can jump non locally on the component in the process of optimization. On such incidences, the gradient based optimizer will immediately reverse the previous geometry change leading to  \emph{ping-pong like}  sequential geometry changes  with essentially no further improvement at all.

In recent years, two of the authors and collaborates have suggested to model low cycle fatigue probabilistically \cite{Gottschalk_Schmitz,Schmitz_Seibel,ASME2013Paper,LCF7Blade,schmitz,bernoulli}, see also \cite{Fedelich,HertelVormwald,Beretta,beretta2,Okeyoyin,amann} for related work. Extensive experimental validation has been provided \cite{Schmitz_Seibel} on various geometries using also notch support factors \cite{maede}. Probabilistic models for low cycle fatigue are natural due to the considerable scatter in LCF life \cite{Vormwald,Harders_Roesler} which often is one order of magnitude. As a byproduct of the probabilistic modelling, the peculiarities of structural integrity with regard to the differentiability in the design parameters is overcome. In this work we will therefore apply the adjoint method to an objective functional and demonstrate the feasibility of the method for 3D trubomachinery components. For related work using different failure mechanisms and simple 2D geometries, see \cite{bolten}. 

The paper is organized as follows: In the following section we shortly recap the probabilistic life calculation from a conceptual and a numerical prospective. In Section \ref{sec:Lagrange} we outline the adjoint method using the Lagrangian approach \cite{nocedal,troeltsch}.  An outline of the numerical implementation strategy and validation work follows in Section \ref{sec:Num}. In Section \ref{sec:Turbo} we demonstrate the viability of the adjoint method for probabilistic LCF for a 3D compressor blisk geometry. We summarize our work and give some outlook to future developments in a conclusion and outlook section.     

%%%%%%%%%%%%%%%%%%%%%%%%%%%%%%%%%%%%%%%%%%%%%%%%%%%%%%%%%%%%%%%%%%%%%%
\section{Computations for Probabilistic Low Cycle Fatigue}
\label{sec:PeobLCF}

This section introduces some notation for probabilistic LCF, see \cite{Vormwald,Harders_Roesler} for an exhaustive treatment. LCF occurs under cyclic loading of technical units. The loading is described by volume forces $f(x)$, like gravity or centrifugal force, surface forces $g(x)$, where $x\in \Omega$ or $x\in\partial\Omega$, respectively. Here $\Omega\subseteq\mathbb{R}^3$ is the region filled with matter. Defining the linearized strain tensor $\varepsilon(u)= \frac{1}{2}(\nabla u+\nabla u^T)$ as the symmetrized Jacobi matrix of the displacement field $u(x)$, $x\in\Omega$ and assuming a linear, isotropic material law $\sigma(u)(x)=\lambda {\rm tr}(\varepsilon)I+2\mu\varepsilon$, where $I$ is the identity matrix and ${\rm tr}$ the trace, we arrive at the usual state equation of linear elasticity in the weak form, see \cite{ciarlet,ern},
\begin{equation}
\label{eqa:Elasticity}
B(u;v):=\int_\Omega \sigma(u):\varepsilon(v)\, dx=F(v):=\int_\Omega f\cdot v\, dx+\int_{\partial\Omega}g\cdot v\, da.
\end{equation}
If \eqref{eqa:Elasticity} is fulfilled for all test functions $v(x)$  (in a suitable function space with three dimensional values), $u(x)$ is a (weak) solution of the elasticity equation. 

To discretize \eqref{eqa:Elasticity}, we choose the node vectors $X$ of a finite element mesh.  $X$ is a tensor or array with dimensions $N\times 3$, where $N$ is the number of (global) nodes. Equivalently, a tensor $X^{\rm loc}$ with the dimensions $M\times q\times 3$ can be used to store the same information, where $M$ is the number of elements and $q$ the order of the finite element. Both formats, as usually, are linked by the connectivity table $l_k(j)$ between local and global degrees of freedom (DoF) $j$ and $l$ given the element $k$.

\begin{figure}[t]
\centerline{\includegraphics[width=.5\textwidth]{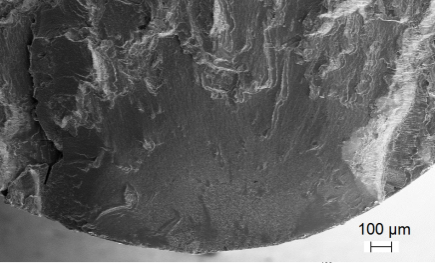}}
\vspace{.2cm}
\caption{Fractographic picture of a LCF crack that originated at the lower surface of the specimen during strain controlled cyclic testing. The rough upper part corresponds to the final violent break after a period of stable crack growth responsible for the lower smooth crack surface. }
\label{fig:LCF}
\end{figure}

Using isoparametric, Lagrangian finite elements  with reference shape functions $\hat\theta_j(\hat x )$ \cite{ern}, we define the mapping between the fixed reference element and the element $k$ in the component via $T_k(\hat x)=\sum_{j=1}^qX^{\rm lok}_{k,j,\cdot}\hat \theta_j(\hat x)$. One defines global shape functions $\theta_l(x)$ for some global node index $l$  such that $\theta_l(x)$ on an element $k$, containing the global node $X_{l,\cdot}$, is given by $\theta_l(x)=\hat\theta_{j}(T^{-1}_k(x))$ where $l=l_k(j)$. In this way, using the finite element approximation $u(x)=\sum_{l=1}^NU_{k,\cdot}\theta_l(x)$ and  \eqref{eqa:Elasticity} can be transformed to the algebraic equation
 \begin{equation}
\label{eqa:StateI}
\sum_{l,r}B_{lr;ms}U_{l,r}=F_{ms}, ~B_{l,r;m,s}=B(\theta_le_r,\theta_m e_s)~~\mbox{and}~~F_{ms}=F(\theta_m e_s),
\end{equation}
with $e_s$ the unit vector with 1 in position $s$ and $0$ elsewhere. The elements of the stiffness matrix $B$ and of the force vector $F$ are now calculated locally, using \eqref{eqa:Elasticity} and numerical quadrature for the volume and surface integrals \cite{ern}.  The result apparently depends on the node coordinates $X$, we write the discretized elasticity equation as 
\begin{equation}
\label{eqa:StateII}
B(X)U(X)=F(X).
\end{equation}

With our notation, we emphasize the dependency of the solution $U=U(X)$ (and other objects introduced below) on the node coordinates $X$. In the following we will call this the \emph{state equation}. Solving the state equation, we obtain the approximated displacement field $u^{\rm fea}(X)(x)=\sum_{l=1}^NU(X)_{l,\cdot}\theta_l(x)$ and the finite element elastic stress field $\sigma^{\rm fea}(X)(x)=\sigma(u^{\rm fea}(X))(x)$. The finite element stress field is now converted to the von Mises equivalent stress amplitude $\sigma_{\rm vM}^{a,\rm fea}(X)(x)$ in the usual way \cite{Harders_Roesler}.

Finishing the stress analysis, the next step is the calculation of the (probabilistic) LCF life. Starting with the deterministic baseline,  $\sigma_{a}^{\rm fea}(X)(x)$ is converted to elastic-plastic equivalent stress using Neuber's rule \cite{Neuber}, see also \cite{Glinka,Knop_Jones_Molent_Wang},
\begin{equation}
\label{eqa:Neuber}
\frac{(\sigma_{a}^{\rm fea})^2}{E}=\frac{(\sigma^{\rm el-pl}_a)^2}{E}+\sigma^{\rm el-pl}_a\left(\frac{\sigma^{\rm el-pl}_a}{K}\right)^{1/n'}.
\end{equation}
Here $E$ is Young's modulus, $K$ is a stress scale for plastic deformation and $n'$ is the hardening exponent \cite{Harders_Roesler}. In \eqref{eqa:Neuber} we suppressed the dependency $\sigma^{\rm el-pl}_a=\sigma^{\rm el-pl}_a(X)(x)$ on the coordinate $x\in\Omega$ and the node set $X$ for notational convenience. We use the notation $\sigma^{\rm el-pl}_a={\rm SD}(\sigma^{\rm fea}_a)$ for the solution of \eqref{eqa:Neuber}.

In the next step we calculate the determinstic LCF cycles to crack initiation, $N_i^{\rm det}(\sigma_a^{\rm el-pl})$. To this aim, we first use the Ramberg-Osgood equation \cite{Ramberg,Harders_Roesler} to convert $\sigma_a^{\rm el-pl}$ to an equivalent strain amplitude
\begin{equation}
\label{eqa:RO}
\varepsilon_a^ {\rm el-pl}={\rm RO}(\sigma_a^{\rm el-pl})=\frac{\sigma_a^{\rm el-pl}}{E}+\left(\frac{\sigma_a^ {\rm el-pl}}{K}\right)^{1/n'}.
\end{equation}  
Next, $N_i^{\rm det}(X)(x)={\rm CMB}^{-1}(\sigma_a^{\rm el-pl}(X)(x))$ is calculated via the Coffin Manson Basquin (CMB) equation \cite{Basquin,Coffin,Harders_Roesler}
\begin{equation}
\label{eqa:CMB}
\varepsilon_a^ {\rm el-pl}={\rm CMB}(N_i^ {\rm det})=\frac{\sigma'_f}{E}(2N_i^ {\rm det})^ b+\varepsilon'_f(2N_i^ {\rm det})^ c.
\end{equation}
$\sigma_f',\varepsilon_f'>0$ and $b,c<0$ are material parameters to be determined from tensile experiments.The deterministic LCF life is now calculated as
\begin{equation}
\label{eqa:DetLife}
N_i^ {\rm det}(X)=\min_{x\in\partial\Omega} N_i^ {\rm det}(X)(x).
\end{equation}
Note that LCF is a surface driven damage mechanism, see Figure \ref{fig:LCF}. Therefore, the minimum over points $x\in\partial\Omega$ is taken. As mentioned in the introduction, considering the weakest spot on the component's surface with the shortest deterministic LCF-life leads to a non differential behavior of $N_i^ {\rm det}(X)$ in the finite element nodes $X$. We also note that the minimum over all $x\in \partial \Omega$ for practical purposes is replaced by the minimum over all nodes $X$ that lie on the surface. 

Let us now contrast this with a probabilistic life calculation as proposed in \cite{Gottschalk_Schmitz,Schmitz_Seibel}. For simplicity, here we do not take into account notch effect modelling \cite{maede}. The probability of failure (PoF) is modelled as a function of the load cycles $n$ 
\begin{equation}
\label{eqa:probLCF}
PoF(X)(n)=1-\exp\left\{-\int_{\partial\Omega}\left(\frac{n}{N_i^{\rm det}(X)(x)}\right)^mda\right\}.
\end{equation}  
Note that the CMB-parameters $\sigma_f'$ and $\varepsilon_f'$ have to be recalibrated when using the probabilistic model \eqref{eqa:probLCF}, see \cite{Schmitz_Seibel}. This results in a Weibull distribution for the probabilistic LCF failure time $N_i^{\rm prob}\sim {\rm Wei}(\eta,m)$, where $m$ is the shape and $\eta=\left[\int_{\partial\Omega}\left(\frac{n}{N_i^{\rm det}(X)(x)}\right)^mda\right]^{-\frac{1}{m}}$ the scale parameter.

 Using finite elements, we now calculate the (approximate) PoF using a numerical quadrature formula for the surface integral
\begin{equation}
J(X,U(X))=\sum_{\mathcal{F}}\sum_{j=1}^{n_q}\hat\omega_j^F\left(\frac{1}{N_i^{\rm det}(X)(T_{k(\cal{ F})}(\hat\xi_j))}\right)^m\sqrt{\hat g_F(X)(\hat\xi_j)},
\end{equation}
where $\hat \omega^F_l$ and $\hat\xi_j$ are the $q^F$ surface quadrature weights on the reference face that corresponds to the face $\mathcal{F}$ on the surface of the component \cite{ASME2013Paper,schmitz}.  $\hat g_F(X)(\hat x)$ is the Gram determinant for the transformation between reference face and the face in the component, which depends on $T_k(X)(\hat x)$ and thereby on the node coordinates $X$. $k=k(\mathcal{F})$ is the element number that contains the surface $\mathcal{F}$. The approximate PoF then is 
\begin{equation}
\label{eqa:probLCFapprox}
PoF(X)(n)\approx 1-\exp\left\{-n^ mJ(X,U(X))\right\}.
\end{equation}
At the end of this section, we have expressed the approximate PoF as a function of the node set $X$ and the finite element global degrees of freedom $U(X)$.

%%%%%%%%%%%%%%%%%%%%%%%%%%%%%%%%%%%%%%%%%%%%%%%%%%%%%%%%%%%%%%%%%%%%%%
\section{Lagrangian Approach to the Adjoint Equation}
\label{sec:Lagrange}
Note that minimization of the PoF \eqref{eqa:probLCFapprox} corresponds to minimization of $J(X,U(X))$ in the discretized geometry $X$. In contrast to the deterministic life $N_i^{\rm det}(X)$, this functional \emph{is} diffenertiable wrt the geometry of the component encoded in the FEA node vectors $X$. As explicit calculations tend to be lengthy, we refer to \cite{saadi} for the details.

At this point, one would like to employ the shape gradient \begin{equation}
\frac{d J(X,U(X))}{dX}=\frac{\partial J(X,U(X))}{\partial X}+\frac{\partial J(X,U(X))}{\partial U}\frac{\partial U(X)}{\partial X}
\end{equation} 
for a gradient based optimization procedure.   However, the computational cost to determine the partial derivatives in $\frac{\partial U(X)}{\partial X}$ is prohibitive, as one FEA calculation would be required for each degree of freedom in $X$. Here we omitted the contraction of various tensor indices for notational simplicity. 

The Lagrangian method helps to circumvent this problem at the cost of one more finite element analysis for what is called the \emph{adjoint equation}. If we consider $U$ and $X$ as two independent sets of variables that enter $J(X,U)$, instead of the minimization of $J(X,U(X))$ in $X$, we consider the minimization of $J(X,U)$ in $U,X$ under the constraint \eqref{eqa:StateII}. We therefore set up the Lagrangian functional \cite{nocedal,troeltsch} using the adjoint state $\Lambda$ -- nothing but the Lagrange multiplier -- as
\begin{equation}
\label{eqa:Lagrange}
\mathcal{L}(X,U,\Lambda)=J(X,U)-\Lambda^ T\left(B(X)U-F(X)\right).
\end{equation}
We now use the Lagrangian formalism. Note that setting the variation of $\mathcal{L}(X,U,\Lambda)$ with respect to $\Lambda$ equal to zero, $\frac{\partial\mathcal{L}}{\partial \Lambda}\stackrel{!}{=}0$ results in the state equation \eqref{eqa:StateII}. Similarly, the adjoint equation is defined by variation of $\mathcal{L}$ wrt the sate variable $U$
\begin{equation}
\label{eqa:Adjoint}
0\stackrel{!}{=}\frac{\partial\mathcal{L}(X,U,\Lambda)}{\partial \Lambda}~~\Leftrightarrow~~B(X)^T\Lambda=\frac{\partial J(X,U)}{\partial U}.
\end{equation} 
The total shape sensitivity  can now be  expressed as the partial derivative of the Lagrangian functional, where $U$ and $\Lambda$ fulfil the state and adjoint equations
\begin{align}
\label{eqa:TotSens}
\begin{split}
\frac{dJ(X,U(X))}{dX}&=\frac{\partial\mathcal{L}(X,U,\Lambda)}{\partial X}\\
&=\frac{\partial J(X,U)}{\partial X}-\Lambda^ T\left(\frac{\partial B(X)}{\partial X}U-\frac{\partial F(X)}{\partial X}\right).
\end{split}
\end{align} 
By \eqref{eqa:probLCFapprox}, the shape sensitivity of the PoF is given by 
\begin{equation}
\frac{dPoF(X)(n)}{dX}=n^m e^{-n^mJ(X,U(X))}\frac{d J(X,U(X))}{dX}.
\end{equation}
The shape sensitivity of $J$ and the PoF thus coincide up to a (global) positive factor. In the follwing we therefore only consider the derivatives of the objective functional $J$. 
\section{Numerical Computation and Validation of Shape Sensitivities }
\label{sec:Num}
 
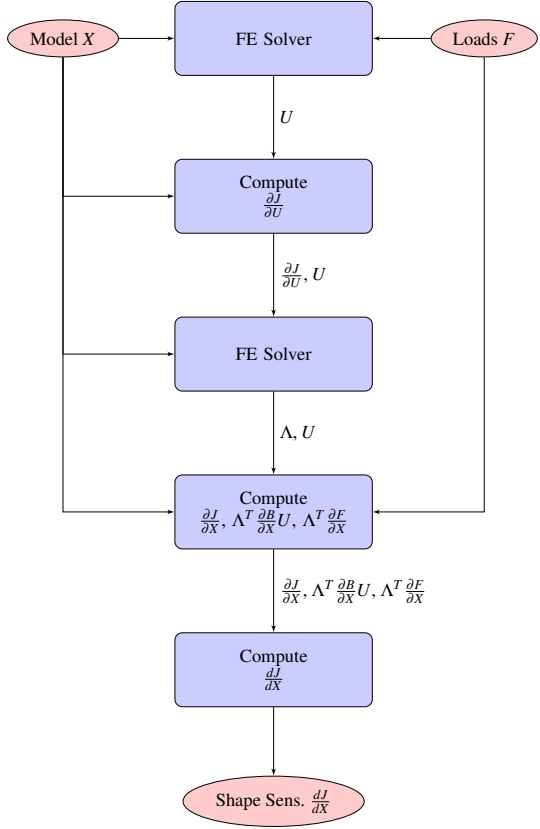
\begin{figure}[t]%\begin{wrapfigure}{rht}[7pt]{0.5\textwidth}
\centerline{\scalebox{0.7}{
\tikzstyle{block} = [rectangle, draw, fill=blue!20, 
    text width=10em, text centered, rounded corners, minimum height=4em, node distance=3cm]
\tikzstyle{line} = [draw, -latex']
\tikzstyle{cloud} = [draw, ellipse,fill=red!20,minimum height=2em,node distance=4cm]
\tikzstyle{cloudL} = [draw, ellipse,fill=red!20,minimum height=2em,node distance=2.5cm]    
\begin{tikzpicture}[node distance = 2cm, auto]
    % Place nodes
    \node [block] (solvU) {FE Solver};
    \node [cloud, left of=solvU] (model) {Model $X$};
    \node [cloud, right of=solvU] (loads) {Loads $F$};
    \node [block, below of=solvU] (cal_dJ_dU) {Compute\\ $\frac{\partial J}{\partial U}$};    
    \node [block, below of=cal_dJ_dU] (solvL) {FE Solver};
	\node [block, below of=solvL] (partial) {Compute\\ $\frac{\partial J}{\partial X}$, $\Lambda^T\frac{\partial B}{\partial X}U$, $\Lambda^T\frac{\partial F}{\partial X}$};  
	\node [block, below of=partial] (final) {Compute\\$\frac{d J}{d X}$};    
	\node [cloudL, below of=final] (output) {Shape Sens.\  $\frac{d J}{d X}$};    
    % Draw edges
    \path [line] (model) --(solvU);
    \path [line] (loads) -- (solvU);
    \path [line] (solvU) -- node  {$U$} (cal_dJ_dU);
	\path [line] (model) |- (cal_dJ_dU);
    \path [line] (model) |- (solvL);
    \path [line] (cal_dJ_dU) -- node  {$\frac{\partial J}{\partial U}$, $U$} (solvL);
    \path [line] (model) |- (partial);
    \path [line] (loads) |- (partial);
    \path [line] (solvL) -- node  {$\Lambda$, $U$} (partial);  
    \path [line] (partial) -- node {$\frac{\partial J}{\partial X}$, $\Lambda^T\frac{\partial B}{\partial X}U$, $\Lambda^T\frac{\partial F}{\partial X}$} (final);
    \path [line] (final) -- (output);
    %\path [line] (loads) |- (solvL);
\end{tikzpicture}
}}
\caption{\label{fig:FowDiag}: Flow chart to calculate the shape sensitivity}
\end{figure}

In this section we describe the numerical implementation and validation for probabilistic LCF. In order to follow the calculations from Setion \ref{sec:Lagrange}, one first has to set up a FEA model, solve the model, extract the node sets and displacements and calculate the objective function $J(X,U)$. the partial $U$ derivatives thereafter have to be calculated and to be fed back into the solver, in order to solve for the adjoint state $\Lambda$.  Thereafter, the quantities $\frac{\partial J}{\partial X}$, $\frac{\partial B}{\partial X}$ and $\frac{\partial F}{\partial X}$ have to be computed and assembled to the shape gradient (or shape sensitivity) \eqref{eqa:TotSens}. Figure \ref{fig:FowDiag} displays this algorithm. 

%\end{wrapfigure}

It is obvious that a naive storage of arrays of dimension $N\times 3\times N\times 3$, like in the case of $\frac{\partial F}{\partial X}$, or even $N\times 3\times N\times 3\times N\times 3$ in the case of $\frac{\partial B}{\partial X}$, exceeds the memory available on usual architectures even for medium FEA-models. Therefore, these objects are calculated for the local degrees of freedom, then are contracted with the local descriptions of $\Lambda^T$ and $U$ element wise and are assembled afterwards.

 \begin{figure*}[ht]
\centerline{\includegraphics[width=.45\textwidth]{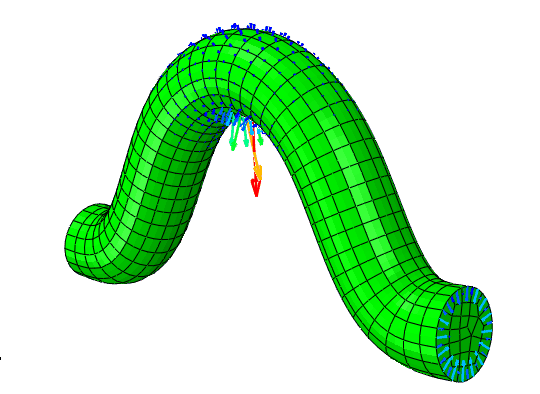}\includegraphics[width=.45\textwidth]{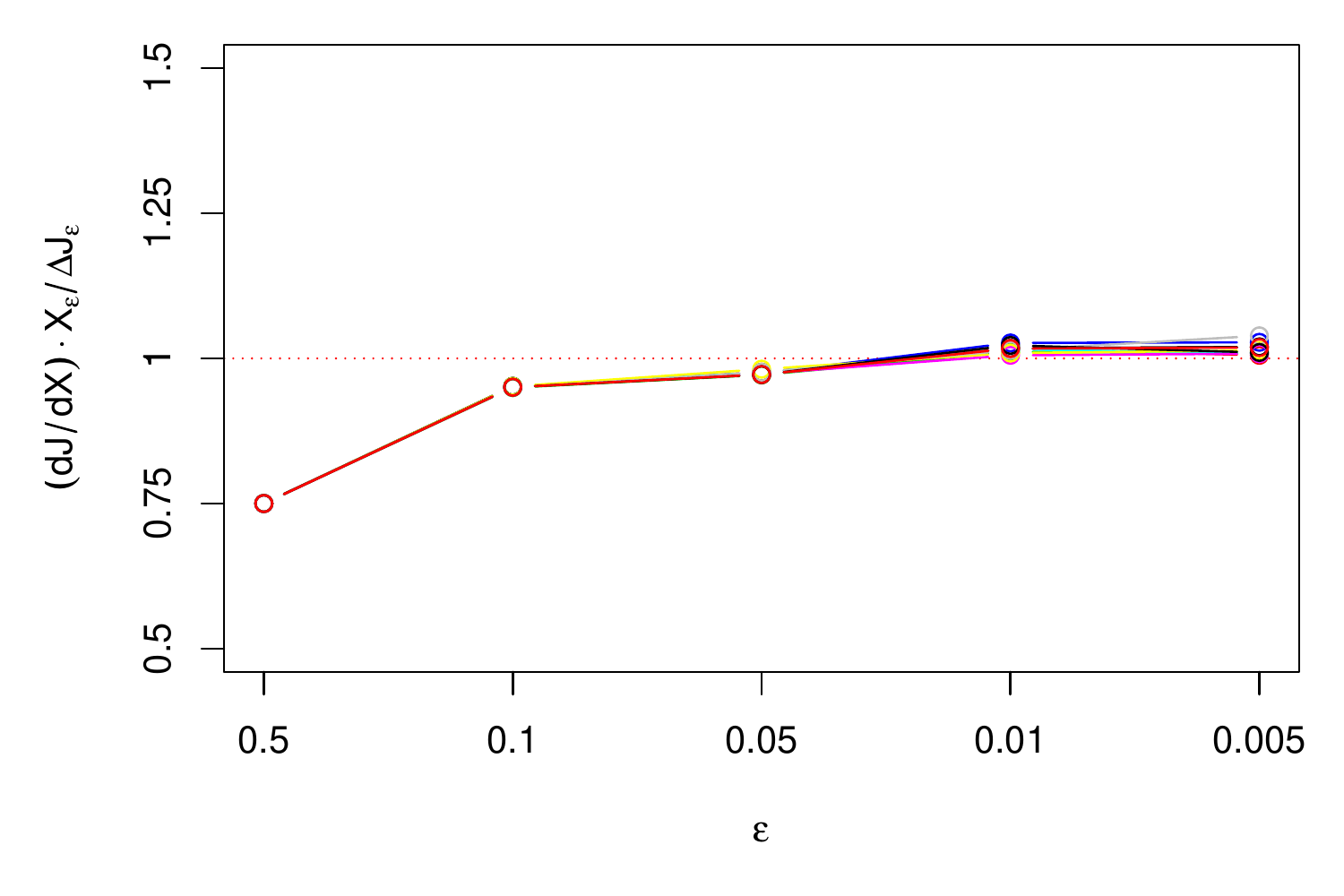}}
\caption{Left panel: Shape sensitivity of LCF failure probability for a bended rod fixed in the rear and pulled on the front face. Right panel: Validadtion results of shape sensitivity divided by finite differences for various stepsizes and random directions. Direction shown is the negative gradient (direction of higher reliability).  }
\label{fig:Rod}
\end{figure*}

We use \texttt{ABAQUS CAE 6.13} as FE solver and  \texttt{R} scripts  (version 3.1.0) for the calculation of partial derivatives and assembly. The element wise local calculations are parallelized with the aid of the \texttt{R}-package \texttt{doParallel}. Interfaces are created between \texttt{ABAQUS} and \texttt{R} to extract model and displacement information and feed back $\frac{\partial J}{\partial U}$ as the right hand side to the FE-solver. This is done on the level of algebraic quantities to avoid the decomposition of the right hand side into surface and volume loads. Dirichlet (\texttt{encastre}) boundary conditions are inherited by the adjoint FE-problem.    

As a validation case, we first set up a 3D model of a bended rod, see Figure \ref{fig:Rod} (left). The geometry is attached to a wall in the rear part and subject to a uniform force density in the front that is pulling the geometry away from the wall. Obviously, the highest stress concentration takes place at the bottom in the middle of the bended region.

The finite element model consists of $N=6410$ nodes distributed over $M=1302$ brick elements of type \texttt{C3D20R} with $q=20$ local degrees of freedom. The reduced quadrature for the assembly of the stiffness matrix $B$ contains $8$ quadrature points and we use $q^F=36$ quadrature points for the surface quadrature \eqref{eqa:probLCFapprox}, see \cite{ASME2013Paper} for a convergence study  wrt the surface quadrature that indicates that refined quadratures are indeed needed. 

The results of the calculation of the total sensitivity are displayed by the arrows in  Figure \ref{fig:Rod}. As one can see, the direction of improved reliability given by the  negative shape gradient points downward and outward. The outward direction aims to diminish the risk of LCF failure by adding more material. The longest arrows in the downward direction in the middle part in the rod clearly aims to reduce stress concentration at the critical spot. See also \cite{bolten} for a 2D counterpart, where actual shape flows under volume constraints have been constructed that ultimately converge to the optimal -- straightened -- configuration.  

Finally, we have to interpret the inward pointing arrows at the surface, where the force is applied. This is explained by the constant force density applied to this surface. Thus, reduction of the surface leads to an effective reduction of the pulling force and thereby to a reduced LCF failure probability. If the force is assumed to be constant, the surface force density $g=g(X)$ has to be readjusted depending on the node set configuration $X$. This leads to extra terms in $\frac{\partial F}{\partial X}$, which change the direction of the shape gradient at the front face of the rod.

Numerical validation work has been conducted by comparison of the shape gradients with finite difference calculations for different stepsizes and five randomly chosen directions of deformation. The right panel in Figure \ref{fig:Rod} reports the results. As we can see, we reach pretty good accuracy $\sim 1\%$ relative error for small stepsizes, where the accuracy for the smallest stepsizes becomes affected by numerical error.

\section{Application to a 3D Turbo Component}
\label{sec:Turbo}

 \begin{figure*}[ht]
\centerline{\includegraphics[width=.45\textwidth]{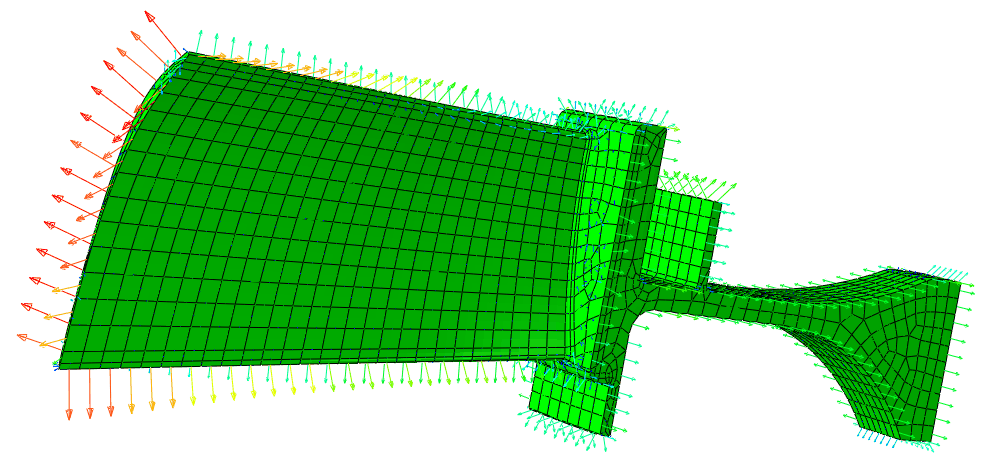}~~~~~\includegraphics[width=.45\textwidth]{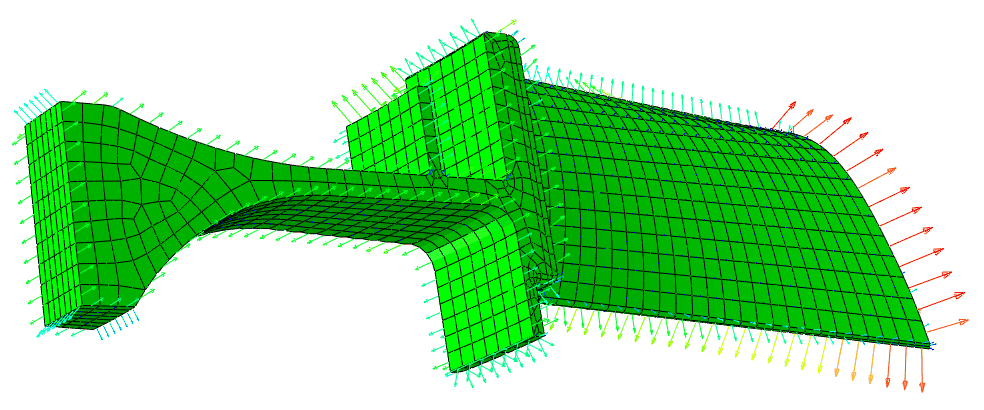}}
\caption{Two views on the shape sensitivity of a 1st stage axial compressor blisk. The arrows show the positive gradient (direction of less reliability). }
\label{fig:Blisk}
\end{figure*}

We finally apply the adjoint method to the shape sensitivity of an aero gas turbine 1st stage compressor blisk model, see Figure \ref{fig:Blisk}. Here we have to account for centrifugal loads and cyclic boundary conditions with an angle $\varphi=2\pi/n_{\rm blades}$, where the blade count is $n_{\rm blade}=45$ in our example. The structure has a height of 38.2 cm and chord length of the blade is 13.88 cm. The bore is 16.9 cm over the rotation axis.  

As a material, we use Titanium Ti-6Al-2Sn-4Zr-2Mo which is regularly applied for high-temperature jet engines. The deterministic CMB parameters are taken from  \cite{heckel}. As explicit experimental data has not been reported, the Weibull shape parameter $m$ is set to a value that matches the usual signal to noise ratio in cyclic tensile testing.       These values are rescaled to probabilistic CMB parameters using the procedure from \cite{Schmitz_Seibel}. The necessary geometric information on the specimens to incorporate the statistical size effect can be found in \cite{heckel2}.

In the case of centrifugal loads, the volume force density is given by
\begin{equation}
\label{eqa:Cerntrifugal}
f(x)=\varrho|x^ \perp|\omega^2,
\end{equation}
where $\varrho$ is the density of the material, $\omega$ the angular velocity and $x^\perp$ is the component of $X$ that is orthogonal to the rotation axis. When calculating shape sensitivities, it has to be taken into account that with a modification of the node set $X$ the volume quadrature points $\xi_V=T_k(X)(\hat\xi_V)$ where $f(\xi_V)$ is evaluated during the assembly of the force vector $F(X)$ change as well. Here $\hat\xi_V$ is the volume quadrature point on the reference element. This creates extra terms in the partial derivatives $\frac{\partial F}{\partial X}$ that have to be properly implemented.

Similarly, the cyclic boundary conditions leas to the indentification of node positions on the front and rear rotor and platform part of the blisk. This has to be taken into account in threefold manner: First, the number of global degrees of freedom is reduced as compared with the unconstrained model. Second, during assembly effects of geometry modification at a face where the cyclic boundary conditions are taken into account have to be trasported from one flank of the structure to the other. When doing this with derivative information, i.e. directional vectors, the proper $2\pi/n_{\rm blades}$ rotations have to be applied. Finally, virtual surfaces on the flanks do not contribute to the probabilistic functional $J(X,U)$. Centrifugal loads are applied at a rotational speed of 397 rad/sec together with \texttt{encastre} boundary conditions on the bore. In the absence of a valid CFD-calculation, we did not apply any pressure loads.   

The blisk model consists of 2262 \texttt{C3D20R} quad elements as in the previous section. The number of global nodes is 13682 and the number of surface faces is 2444. The following Table \ref{tab:Runtime} reports execution times for the critical steps of Figure \ref{fig:FowDiag} for the present model.
\begin{table}[h]
\centerline{
\begin{tabular}{|l|c|c|l|}
\hline
\textbf{Quantity}&\textbf{Elapsed}&\textbf{Cores}&\textbf{Tool}\\\hline
State $U$&31.52&1&\texttt{ABAQUS CAE 6.13}\\
$\frac{\partial J}{\partial U}$,$\frac{\partial J}{\partial X}$ together&14.38&1&\texttt{R 3.1.0}\\
Adjoint State $\Lambda$&37.49&1&\texttt{ABAQUS CAE 6.13}\\
$\frac{\partial B}{\partial X}$,$\frac{\partial F}{\partial X}$ together&378.52&6&\texttt{R 3.1.0}\\\hline
\end{tabular}
\vspace{.2cm}
}
\caption{Execution times in sec on an  Intel Core i7-3630QM CPU @ 2.40GHZ, 8GB shared memory machine with 4 physical and 8 virtual cores.} 
\label{tab:Runtime}   
\end{table}

We see that the main time consumption is needed to calculate the last term in \eqref{eqa:TotSens} which consumed more than 6 minutes for our software prototype. The reason is that even single element shape derivatives of the stiffness matrix require the calculation of an array of dimension $q\times3\times q\times 3\times q\times 3$ times the number of volume quadrature points, which for an element with 20 DoF and 8 volume quadrature points for a reduced quadrature  requires the calculation of more than 1.728 million array entries per element. Element types with less degrees of freedom however lead to much more noisy representations of the shape sensitivities.      

We note that these leading local computations parallelize without significant overhead and can thus be reduced proportionally to the number of cores. Secondly, the complexity of the calculation of these terms scales linearly in the number of elements and will thus become sub-leading for models of larger size as in our example. Also, portation to compiled code has a certain potential, although the workhorse, \texttt{R}'s array arithmetic, is compiled \texttt{C} and \texttt{FORTRAN} code.    

Let us now shift attention to physical interpretation of the calculation procedure. Figure \ref{fig:Blisk} clearly shows by the size of the outward pointing arrows (direction of less reliability) that adding more material to the rotating system at almost all locations will lead to less reliability due to higher centrifugal loads. Consistently, these effects get worse as the length of the arrows increase the more one approaches the blade tip. In other locations, as the fillets below the platform, we observed inward pointing gradients suggesting that more material in this region, despite higher centrifugal forces, is capable to improve the design's reliability.

Let us also note that the outward arrows at the flanks identified via the circular boundary conditions are of artificial nature as the discretized model only respects these conditions via inter-nodal constraints. Therefore, the model artificially predicts shape sensitivities at these flanks which of course do not correspond to any design option. This artefact can be resolved by simply setting these sensitivities to zero.

%%%%%%%%%%%%%%%%%%%%%%%%%%%%%%%%%%%%%%%%%%%%%%%%%%%%%%%%%%%%%%%%%%%%%%
\section{Conclusion and Outlook}
\label{sec:Conclusion}
In this work we have shown that a probabilistic description of probabilistic low cycle fatigue has the beneficial side effect of enabling the adjoint method for the calculation of shape sensitivities. This method has been validated for 3D examples and applied to realistic 3D turbo geometries. 

While we have demonstrated that our approach is viable, further development work has to be conducted in order to fully exploit the potential of the method for gas turbine design. 

First, the probabilistic functional should be extended to contain also notch effects \cite{maede}. While in principle there is no problem to replace $J(X,U)$ with a more accurate form, numerical difficulties might arise from the use of second order derivatives of $u(x)$ in the calculation of the local notch support  $\chi^*(x)$ on the component's surface.
Secondly, the mechanical finite element model has to be extended to a thermomechanical model in order to deal also with cooled turbine components \cite{Hetnarski,bittner}. Also, the treatment of contact boundary conditions still has to be integrated \cite{wehrstedt}, with some impact on the Lagrangian formalism due to inequality constraints instead of the state equation constraint \eqref{eqa:StateII} of equality type. 
Furthermore, the backreaction of the flow field to the change of geometry will also change the surface pressure $g(X)$, such that fluid structure interaction has to be taken into account.  

Finally, to return to our initial motivation, the minimization of the PoF alone may lead to shapes that have a poor performance in other relevant criteria like, for example, efficiency, volume and cost.

\begin{figure}[ht]
\begin{center}
\includegraphics[scale=0.54]{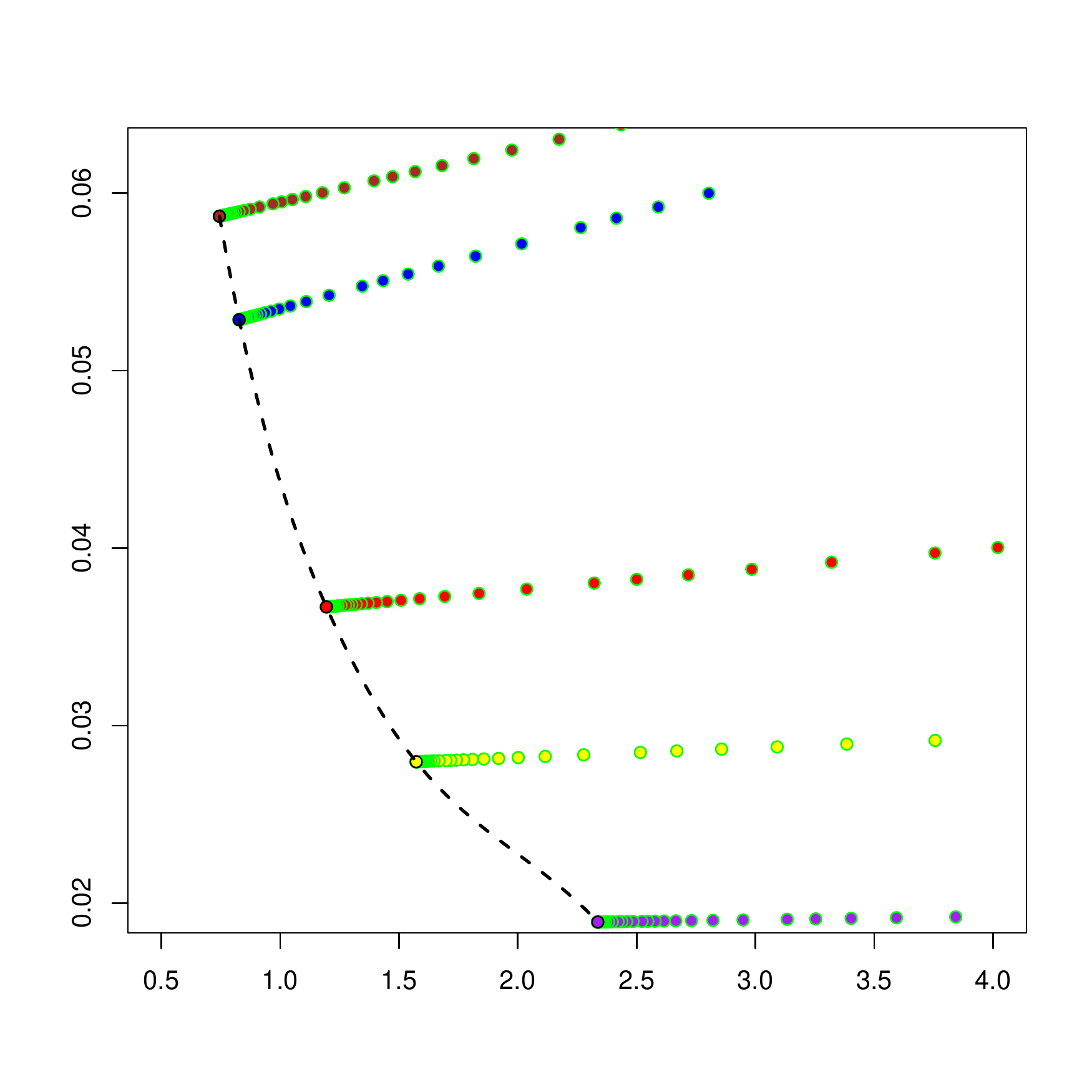}
\end{center}
\caption{Exploration of the Pareto front using a gradient based multicriteria steepest decent algorithm \cite{fliege} for a 2D rod geometry with objective functions probability of failure (PoF) and material consumption (volume), see \cite{doganay}.  }
\label{Pareto-Front}
\end{figure}

The availability of shape gradients for turbomachinery components facilitates the simultaneous consideration of several optimization criteria. 
As a preliminary two-dimensional case study, we considered the simultaneous minimization of the PoF \emph{and} the volume of a 2D ceramic rod (similar to Fig. \ref{fig:Rod}) in a biobjective shape optimization problem. Since in this context there usually does not exist a shape minimizing both objective functions at the same time, we aim at the computation of \emph{Pareto optimal} shapes that can not be improved in one objective without deterioration in the other objective.
In other words, if $X^*$ is a Pareto optimal discretized geometry, then there exists no other discretized geometry that is in all objective functions better or equal than $X^*$, and in at least one objective function strictly better.

We have implemented a multi-objective descent algorithm \cite{fliege} that can be interpreted as a generalization of a classical steepest descent method to the case of multiple objective functions. Given the shape gradients of the two objective functions PoF and volume, a joint direction of ``steepest descent'' is computed in each iteration using an auxiliary quadratic optimization problem. %The steplength is computed with an Armijo-like rule, extended by regularity constraints.

Preliminary computational tests show a fast convergence towards the Pareto front. An approximation is obtained by varying the starting shape, see Figure \ref{Pareto-Front}. This approximation provides information on the trade-off between the two objective functions, supporting the decision making process and suggesting several interesting solution alternatives. An extension to three-dimensional shapes and further objective functions such as efficiency is currently under research.

%%%%%%%%%%%%%%%%%%%%%%%%%%%%%%%%%%%%%%%%%%%%%%%%%%%%%%%%%%%%%%%%%%%%%%
\begin{acknowledgment}
We thank the Siemens Gas Turbine Technology Department, in particular Dr. Georg Rollmann, for constant support. Also we thank Matthias Bolten, Camilla Hahn, Nadine Moch, Lucas M\"ade, Benedict Engel and Tilmann Beck for illuminating discussions on probabilistic LCF and the permission to use Figure \ref{fig:LCF}. We gratefully acknowledge the permission to use the blisk geometry shown in Figure \ref{fig:Blisk} designed by BWEngineering in this work.   This article has grown out of a AG Turbo project sponsored by the federal ministry of economic affairs (BMWi, grant-no 03 ET 2013I) and the federal ministry of research and education (BMBF, grant-no: 05M18PXA ) as a part of the GIVEN collaboration.   
\end{acknowledgment}

%%%%%%%%%%%%%%%%%%%%%%%%%%%%%%%%%%%%%%%%%%%%%%%%%%%%%%%%%%%%%%%%%%%%%%
% The bibliography is stored in an external database file
% in the BibTeX format (file_name.bib).  The bibliography is
% created by the following command and it will appear in this
% position in the document. You may, of course, create your
% own bibliography by using thebibliography environment as in
\bibliographystyle{asmems4}

%\vspace{45cm}
%${}$
%\vspace{18cm}

\noindent \textbf{Permission for Use} - The content of this paper is copyrighted by Siemens Energy, Inc. and is licensed to ASME for publication and distribution only. Any inquiries regarding permission to use the content of this paper, in whole or in part, for any purpose must be addressed to Siemens Energy, Inc. directly. 

\noindent The authors are granted the right to self-archive.

% Here's where you specify the bibliography database file.
% The full file name of the bibliography database for this
% article is asme2e.bib. The name for your database is up
% to you.
%\bibliography{asme2e}

\begin{thebibliography}{99}
\bibitem{amtsfeld} P. Amtsfeld, M. Lockan, D. Bestle and M. Meyer , Accelerated 3D Aerodynamic Optimization of Gas Turbine Blades, ASME Turbo Expo 2014,  GT2014-25618.
\bibitem{kim} Y. Kim, S. Lee, K. Yee, Y.-S. Kang, Aerodynamic Efficiency Optimization of the 1st Stage of Transonic High Pressure Turbine through Lean and Sweep Angles, Int. J. of Turbo \& Jet-Engines, to appear.
 \bibitem{buske}  C. Buske, A. Krumme, T. Schmidt, C. Dresbach, S. Zur and R. Tiefers , Distributed Multidisciplinary Optimization of a Turbine Blade Regarding Performance, Reliability and Castability , ASME Turbo Expo 2017,  GT2016-56079.
\bibitem{backhaus} J. Backhaus, M. Aulich,  C. Frey, T. Lengyel and C. Vo\ss
Gradient Enhanced Surrogate Models Based on Adjoint CFD Methods for the Design of a
Counter Rotating Turbofan. ASME Turbo Expo 2012
\bibitem{xu} G. Yu and F. Christakupoulos,  CAD-Based Adjoint Shape Optimisation of a One-Stage Turbine With Geometric Constraints, ASME Turbo Expo 2015, GT2015-42237
\bibitem{forrester} A. Forrester, A. Sobester and A. Keane, Engineering Design via Surrogate Modeling: a Practical Guide, Wiley 2008.
\bibitem{debar} J.H.S. de Baar; T.P. Scholcz; R.P. Dwight,  Exploiting Adjoint Derivatives in High-Dimensional Metamodels". AIAA Journal. 53 (5) (2015): 1391-1395.
\bibitem{nocedal} J. Nocedal, S. Wright, Numerical Optimization, Springer 2006.
\bibitem{troeltsch} F. Tr\"oltsch, Optimal Control of Partial Differentiel Equations, AMS Graduate Texts in Mathematics, AMS 2010.
\bibitem{giles1} M.B. Giles and N.A. Pierce `Adjoint Equations in CFD: Duality, Boundary Conditions and Solution Behaviour'. AIAA Paper 97-1850, 1997.
\bibitem{giles2} M.B. Giles, M.C. Duta, J.-D. Muller and N.A. Pierce, Algorithm Developments for Discrete Adjoint Methods, AIAA Journal, 41(2), 2003.
\bibitem{wever} S. K\"ammerer, J. Mayer, H. Stetter, M. Paffrath, U. Wever and A. R. Jung, Development of a Three Dimensional Geometry Optimization Methods for Turbomachinery Applications, Int. J. of Rotating Machinery 10 (5) 2004, 373-385.
\bibitem{frey} C. Frey, D. N\"urnberger, and H.P. Kersken. The Discrete Adjoint of a Turbo- Machinery RANS Solver, ASME-GT2009, 2009.
\bibitem{schulz} V. Schulz, M. Siebenborn, Computational Comparison of Surface Metrics
for PDE Constrained Shape Optimization, Computational Methods in Applied Mathematics 16 (3) 2016.
  \bibitem{Gottschalk_Schmitz} H. Gottschalk and S. Schmitz, Optimal Reliability in Design for Fatigue Life, Part I -- Existence of Optimal Shapes, SIAM J. Control Optim., 52(5), pp. 2727--2752, 2014
  \bibitem{Schmitz_Seibel} S. Schmitz, T. Seibel, T. Beck, R. Rollmann, R. Krause and H. Gottschalk, A Probabilistic Model For LCF, Comp. Materials Science 79, 2013, 584--590.
  \bibitem{ASME2013Paper} S. Schmitz, H. Gottschalk, R. Rollmann and R. Krause, 2013, Risk Estimation for LCF Crack Initiation, ASME Turbo Expo  GT2013-94899.
  \bibitem{LCF7Blade} S. Schmitz, R. Rollmann, H. Gottschalk and R. Krause, Probabilistic Analysis of LCF Crack Initiation Life of a Turbine Blade under Thermomechanical Loading, Proc. Int. Conf LCF 7, 2013.
  \bibitem{schmitz} S. Schimitz, A Local and Probabilistic Model for Low-Cycle Fatigue.
New Aspects of Structural Analysis, Hartung-Gorre Verlag, 2015.
\bibitem{bernoulli} L. Bittner , H. Gottschalk, M. Gröger,
N. Moch, M. Saadi and S. Schmitz, Modeling, Minimizing and Managing the Risk of
Fatigue for Mechanical Components, in: S. Albeverio, D. Holms and A. Cruceiro (Ed.) Stochastic Geometric Mechanics - A Series of Lectures, Springer 2017.
  	 \bibitem{Fedelich} B. Fedelich, A Stochastic Theory for the Problem of Multiple Surface Crack Coalescence, Int. J. of Fracture 91,1998,  23--45.
\bibitem{HertelVormwald} O. Hertel and M. Vormwald, Statistical and Geometrical Size Effects in Notched Members Based on Weakest-Link and Short-Crack Modelling, Engineering Fracture Mechanics 95, 2012, 72--83.
	\bibitem{Beretta} S. Beretta, H. J. Desimone and A. Poli, Fatigue Assesment of Tubular Automotive Components in Presence of Inhomogenities, Proceedings of IMECE2004-60333, 2004, 791--798.
	\bibitem{beretta2} S.P. Zhu, S. Foletti and S.Beretta, Probabilistic
Framweork for Multiaxial LCF Assessment under Material Variability,
International Journal of Fatigue 2017 (in press).
	 \bibitem{Okeyoyin} O.A. Okeyoyin, G.M. Owolabi, Application of Weakest Link Probabilistic Framework for Fatigue Notch Factor to Turbine Engine Materials, 13th International Conference on Fracture, 2013.
		\bibitem{amann} C. Amann, K. Kadau,  Numerically Efficient Modified Runge-Kutta
Solver for Fatigue Crack Growth Analysis, Engineering Fracture Mechanics, 161 ,2016., 55–62
 \bibitem{maede} L. Maede, S. Schmitz, H. Gottschalk and T. Beck, Combined Notch and Size Effect Modeling in
a Local Probabilistic Approach for LCF, Comp. Materials Science (2018), to appear.
\bibitem{Vormwald} D. Radaj and M. Vormwald,  Fatigue Resistance (in German), 3rd edition, Springer Berlin Heidelberg, 2007.
  \bibitem{Harders_Roesler} M. B\"aker, H. Harders and J. R\"osler,  Mechanical Behaviour of Engineering Materials: Metals, Ceramics, Polymers and Composites, 1st edition, Springer Berlin Heidelberg New York 2007.
\bibitem{bolten} M. Bolten, H. Gottschalk, C. Hahn and M. Saadi, Shape Optimization to Decrease Failure Probabilty, Preprint 2017, arXiv:1705.05776.
\bibitem{ciarlet} P. Ciarlet, Mathematical Elasticity - Volume I: Three-Dimensional Elasticity, North-Holland, Amsterdam, 1988
    \bibitem{ern} A. Ern and J.-L. Guermond, Theory and Practice of Finite Elements, Springer, New York, 2004.
   \bibitem{Ramberg} W. Ramberg and W. R. Osgood, Description of Stress-Strain Curves by Three Parameters, Technical Notes - National Advisory Committee For Aeronautics, No. 902, Washington DC., 1943
  \bibitem{Neuber} H. Neuber, Theory of Stress Concentration for Shear-Strained Prismatical Bodies with Arbitrary Nonlinear Stress-Strain Law, J. Appl. Mech. 26, 544, 1961.
  \bibitem{Glinka} G. Glinka, Energy Density Approach to Calculation of Inelastic Stress-Strain Near Notches and Cracks, Engineering Fracture Mechanics 22(3), 1985 485--508, .
	\bibitem{Knop_Jones_Molent_Wang} M. Knop, R. Jones, L. Molent, L. Wang, On Glinka and Neuber Methods for Calculating Notch Tip Strains under Cyclic Load Spectra, Int. J. of Fatigue, Vol. 22,  2000,  743--755.
		\bibitem{Basquin} O.H. Basquin, The exponential law of endurance tests, Proc. ASTM, 10, 1910, 625--630.
	\bibitem{Coffin} J. Coffin L. F., A Study of the Effects of Cyclic Thermal Stresses on a Ductile Metal, Trans. ASME 76, 1954, 931--950.
	\bibitem{saadi} M. Saadi, PhD Thesis in Mathematics, University of Wuppertal 2018, in preparation
	\bibitem{heckel} T. K. Heckel, H.-J. Christ, Isothermal and Thermomechanical Fatigue of Titanium Alloys, Procedia
Engineering 2, 2010, 845--854. 
\bibitem{heckel2} T. Heckel, Isothermal and Thermomechanical Fatigue of Titanium Alloys (in German). Doctorate Thesis, University of Siegen. 
Shaker-Verlag 2010.
    \bibitem{Hetnarski} R. B. Hetnarski and M. Reza Eslami, Thermal Stresses - Advanced Theory and Applications, Springer, Berlin, 2009.
    \bibitem{bittner} L. Bittner and H. Gottschalk, Optimal Reliability for Components under Thermomechanical Cyclic Loading, Control \& Cybernetics 45 2016, 2--35.
    \bibitem{wehrstedt} J. C. Wehrstedt, Shape Optimzation with Variatioal Inequalities as Constraint and an Application in Pine Surgery (in German), PhD Thesis in Mathematics, Technical University M\"unchen 2007.
    \bibitem{fliege} J. Fliege, B. F. Svaiter, Steepest Descent Method for Multicriteria Optimization, Math. Methods of Operations Research
51 (3), 2000, 479--494.
\bibitem{doganay} O. T. Doganay, Multicriteria Optimization with Shape
Gradients, Master Thesis Wuppertal 2017.
\end{thebibliography}

%%%%%%%%%%%%%%%%%%%%%%%%%%%%%%%%%%%%%%%%%%%%%%%%%%%%%%%%%%%%%%%%%%%%%%
%\appendix       %%% starting appendix
%\section*{Appendix A: Head of First Appendix}
%Avoid Appendices if possible.
%%%%%%%%%%%%%%%%%%%%%%%%%%%%%%%%%%%%%%%%%%%%%%%%%%%%%%%%%%%%%%%%%%%%%%%
%\section*{Appendix B: Head of Second Appendix}
%\subsection*{Subsection head in appendix}
%The equation counter is not reset in an appendix and the numbers will
%follow one continual sequence from the beginning of the article to the very end as shown in the following example.
%\begin{equation}
%a = b + c.
%\end{equation}

\end{document}